\documentclass[reqno,oneside,12pt]{amsart}
\newcommand{\K}{\overline{K}}
\begin{document}
\vspace*{8pc}
\begin{flushleft}
{\LARGE \bf
Fractional Differentiation Operator over an Infinite Extension of
a Local Field }
\end{flushleft}

\vspace*{2pc}
\noindent
ANATOLY N. KOCHUBEI \quad
Institute of Mathematics,Ukrainian National Academy of Sciences,
Tereshchenkivska 3, Kiev, 252601 Ukraine

\vspace*{5pc}
\noindent {\bf 1\quad INTRODUCTION}

\vspace{2pc}
\noindent
Let $k$ be a non-Archimedean local field of zero characteristic.
Consider an increasing sequence of its finite extensions
$$
k=K_1\subset K_2\subset \ldots \subset K_n\subset \ldots \ .
$$
The infinite extension
$$
K=\bigcup \limits _{n=1}^\infty K_n
$$
may be considered as a topological vector space over $k$ with the
inductive limit topology. Its strong dual $\overline{K}$ is the
basic object of the non-Archimedean infinite-dimensional analysis
initiated by the author [1]. Let us recall its main
constructions and results.

Consider, for each $n$, a mapping $T_n:\ K\to K_n$ defined as follows.
If $x\in K_\nu $, $\nu >n$, put
$$
T_n(x)=\frac{m_n}{m_\nu }\mbox{Tr}_{K_\nu /K_n}(x)
$$
where $m_n$ is the degree of the extension $K_n/k$,
$\mbox{Tr}_{K_\nu /K_n}:\ K_\nu \to K_n$ is the trace mapping. If
$x\in K_n$ then, by definition, $T_n(x)=x$. The mapping $T_n$ is
well-defined and $T_n\circ T_\nu =T_n$ for $\nu >n$. Below we shall
often write $T$ instead of $T_1$.

The strong dual space $\K $ can be identified with the projective
limit of the sequence $\{ K_n\} $ with respect to the mappings
$\{ T_n\}$, that is with the subset of the direct product $\prod
\limits _{n=1}^\infty K_n$ consisting of those $x=(x_1,\ldots
,x_n,\ldots )$,  $x_n\in K_n$, for which $x_n=T_n(x_\nu )$ if
$\nu >n$. The topology in $\K$ is introduced by seminorms
$$
\Vert x\Vert _n=\Vert x_n\Vert ,\quad n=1,2,\ldots ,
$$
where $\Vert \cdot \Vert $ is the extension onto $K$ of the
absolute value $|\cdot |_1$ defined on $k$.

The pairing between $K$ and $\K$ is defined as
$$
<x,y>=T(xy_n)
$$
where $x\in K_n\subset K,\ y=(y_1,\ldots ,y_n,\ldots )\in \K ,\
y_n\in K_n$. Both spaces are separable, complete, and reflexive.
Identifying an element $x\in K$ with $(x_1,\ldots ,x_n,\ldots
)\in \K$ where $x_n=T_n(x)$, we can view $K$ as a dense subset of
$\K$. The mappings $T_n$ can be extended to linear continuous
mappings from $\K$ to $K_n$, by setting $T_n(x)=x_n$ for any
$x=(x_1,\ldots ,x_n,\ldots )\in \K$.

Let us consider a function on $K$ of the form
$$
\Omega (x)=\left\{ \begin{array}{rl}
1 ,& \mbox{if } \|x\|\le 1;\\
0, & \mbox{if } \|x\|>1.
\end{array}\right.
$$
$\Omega$ is continuous and positive definite on $K$. That results in the
existence of a probability measure $\mu $ on the Borel $\sigma$-algebra
$\mathcal B(\K)$, such that
$$
\Omega (a)=\int \limits _{\K}\chi (<a,x>)\,d\mu (x),\quad a\in K,
$$
where $\chi $ is a rank zero additive character on $k$. The measure $\mu $
is Gaussian in the sense of Evans [2]. It is concentrated on the
compact additive subgroup
$$
S=\left\{ x\in \overline{K}:\ \ \|T_n(x)\|\le q_n^{d_n/m_n}\|m_n\|,\
n=1,2,\ldots \right\}
$$
where $q_n$ is the residue field cardinality for the field $K_n$,
$d_n$ is the exponent of the different for the extension $K_n/k$.
The restriction of $\mu $ to $S$ coincides with the normalized Haar
measure on $S$. On the other hand, $\mu $ is singular with respect to
additive shifts by elements from $\K \setminus S$.

Having the measure $\mu $, we can define a Fourier transform
$\widehat{f}=\mathcal Ff$ of a complex-valued function $f\in L_1(\K)$
as
$$
\widehat{f}(\xi )=\int _{\overline{K}}\chi (\langle \xi ,x\rangle )
f(x)\,d\mu (x).
$$

Let $\mathcal E(\K)\subset L_1(\K)$ be the set of ``cylindrical''
functions of the form $f(x)=\varphi (T_n(x))$ where $n\ge 1$,
$\varphi $ is a locally constant function. The fractional
differentiation operator $D^\alpha $, $\alpha >0$, is defined on
$\mathcal E(\K)$ as $D^\alpha =\mathcal F^{-1}\Delta ^\alpha
\mathcal F$ where $\Delta ^\alpha $ is the operator of
multiplication by the function
$$
\Delta ^\alpha (\xi )=\left\{ \begin{array}{rl}
\|\xi \|^\alpha , & \mbox{if }\|\xi \|>1\\
0, & \mbox{if }\|\xi \|\le 1
\end{array} \right. ,\ \ \xi \in K.
$$
Correctness of this definition follows from the theorem of
Paley-Wiener-Schwartz type for the transform $\mathcal F$.
The operator $D^\alpha $ is essentially self-adjoint on
$L_2(\K)$.

The operator $D^\alpha $ is an infinite-dimensional counterpart of the
$p$-adic fractional differentiation operator introduced by
Vladimirov [3] and studied extensively in [4-14]. In some
respects the infinite-dimensional $D^\alpha $ resembles its
finite-dimensional analogue though it possesses some new
features. For example, we shall show below that the structure of
its spectrum depends on arithmetic properties of the extension
$K$.

Both in the finite-dimensional and infinite-dimensional cases
$D^\alpha $ admits a probabilistic interpretation. Namely, $-D^\alpha $
is a generator of a cadlag Markov process, which is a $p$-adic
counterpart of the symmetric stable process. In analytic terms,
that corresponds to a hyper-singular integral representation of
$D^\alpha f$, for a suitable class of functions $f$. In [1]
 such a representation was obtained for $f\in
\mathcal E(\K)$, $f(x)=\varphi (T_n(x))$, $\varphi $ locally
constant:
$$
(D^\alpha f)(x)=\psi (T_n(x))
$$
where
\begin{multline}
\psi (z)=q_n^{d_n\alpha /m_n}\frac{1-q_n^{\alpha /m_n}}{1-q_n^{-1-
\alpha /m_n}}\|m_n\|^{-m_n}\\ \times \int \limits_{x\in K_n:\  \|x\|
\le q_n^{d_n/m_n}\|m_n\|}\left[ \|x\|^{-m_n-\alpha }\|m_n\|^{m_n+
\alpha } +\frac{1-q_n^{-1}}{q_n^{\alpha /m_n}-1}q_n^{-d_n(1+\alpha /m_n)}
\right] \\ \times [\varphi (z-x)-\varphi (z)]\,dx,
\end{multline}
$z\in K_n,\ \ \|z\|\le q_n^{d_n/m_n}\|m_n\|$.

In this paper we shall show that there exists also a
representation in terms of the function $f$ itself:
\begin{equation}
(D^\alpha f)(y)=\int \limits _{\K}[f(y)-f(x+y)]\Pi (dx)
\end{equation}
where $\Pi $ is a measure on $\mathcal B(\K \setminus \{ 0\} )$
finite outside any neighbourhood of zero.

Another measure of interest in this context is the heat measure
$\pi (t,dx)$ corresponding to the semigroup $\exp (-tD^\alpha )$,
$t>0$. We show that in contrast both to the finite-dimensional
case and to the similar problem for the real infinite-dimensional
torus [15], $\pi (t,dx)$ is not absolutely continuous
with respect to $\mu $, whatever is the sequence $\{ K_n\}$.

\vspace{3pc}
\noindent {\bf 2\quad SPECTRUM}

\vspace{2pc}
\noindent
We shall preserve the notation $D^\alpha $ for its closure in
$L_2(\K)=L_2(S,d\mu )$. It is clear from the definition that the
spectrum of $D^\alpha $ coincides with the closure in $\mathbb R$
of the range of the function $\Delta ^\alpha (\xi )$, $\xi \in
K$. In order to investigate the structure of the spectrum, we
need some auxiliary results.

It follows from the duality theory for direct and inverse limits
of locally compact groups [16] that the character group
$\K ^*$ of the additive group of $\K$ is isomorphic to $K$. The
isomorphism is given by the relation
\begin{equation}
\psi (y)=\chi (<a_\psi ,y>),\quad y\in \K,
\end{equation}
where $a_\psi \in K$ is an element corresponding to the character
$\psi $.

Denote $O=\{ \xi \in K:\ \Vert \xi \le 1\} ,\ \ O_n=O\cap K_n$.

\vspace{2pc}
\noindent LEMMA 1\quad
The dual group $S^*$ of the subgroup $S\subset \K$ is isomorphic
to the quotient group $K/O$.

\vspace{1pc}
\noindent Proof: \quad
We may assume (without restricting generality) that $k=Q_p$. Let
$$
S^{(n)}=\{ z\in K_n:\ \Vert z\Vert \le q_n^{d_n/m_n}\Vert
m_n\Vert \} ,\quad n=1,2,\ldots .
$$
It is clear that $S^{(n)}$ is a compact (additive) group. If
$z\in S^{(\nu )}$, $\nu >n$, then
$$
\Vert T_n(z)\Vert =\Vert m_n\Vert \cdot \Vert m_\nu \Vert
^{-1}\Vert \mbox{Tr}_{K_\nu /K_n}(z)\Vert
$$
where $|m_\nu ^{-1}z|_\nu \le q_\nu ^{d_\nu }$, and $|\cdot |_\nu
$ is the normalized absolute value on $K_\nu $ [17]. Therefore (see
Chapter 8 in [18])
$$
|m_\nu ^{-1}\mbox{Tr}_{K_\nu /K_n}(z)|_n\le q_n^l
$$
where $l\in \mathbb Z,\ e_{n\nu }(l-1)<d_\nu -d_{n\nu }\le e_{n\nu }l,\
e_{n\nu }$ and $d_{n\nu }$ are the ramification index and the exponent
of the different for the extension $K_\nu /K_n$. On the other
hand, $d_\nu =e_{n\nu }d_n+d_{n\nu }$, so that $l=d_n$, and we
find that $T_n(z)\in S^{(n)}$.

Hence, $T_n:\ S^{(\nu )}\to S^{(n)}$. It is clear that $T_n$ is a
continuous homomorphism. As a result, $S$ can be identified with
an inverse limit of compact groups $S^{(n)}$ with respect to the
sequence of homomorphisms $T_n$. Using the auto-duality of each
field $K_n$, we can identify the group dual to $S^{(n)}$ with
$K_n/\Phi _n$ where $\Phi _n$ is the annihilator of $S^{(n)}$ in
$K_n$. On the other hand, $\Phi _n=O_n$.

Indeed,
$$
\Phi _n=\{ \xi \in K_n:\ \chi (T(z\xi ))=1\ \forall z\in S^{(n)}\}
$$
$$
=\{ \xi \in K_n:\ |T(z\xi )|_1\le 1\ \forall z\in S^{(n)}\}.
$$
If $\xi \in O_n$ then $|m_n^{-1}z\xi |_n\le q_n^{d_n}$ for any
$z \in S^{(n)}$, whence $\xi \in \Phi _n$ by the definition of
the number $d_n$.

Thus $O_n\subset \Phi _n$. Conversely, suppose that $\xi \in \Phi
_n\setminus O_n$, that is $\Vert \xi \Vert >1$, $|T(z\xi )|_1\le
1$ for all $z\in S^{(n)}$. Let $z$ be such that $\Vert z\Vert
=\Vert m_n\Vert q_n^{d_n/m_n}$. Then
$$
|m_n^{-1}|_n\cdot |\xi |_n\cdot |z|_n=|\xi
|_nq_n^{d_n}>q_n^{d_n}.
$$
It follows from the properties of trace [18] that $z$ can
be chosen in such a way that $|\mbox{Tr}_{K_n/k}(m_n^{-1}\xi
z)|_1>1$, and we have a contradiction. So $\Phi _n=O_n$.

It follows from the identity
$$
T(\xi T_n(z))=T(\xi z),\quad \xi \in K_n,\ z\in K_\nu ,\ \nu >n,
$$
that the natural imbeddings $K_n/O_n\to K_\nu /O_\nu$, $\nu >n$, are the
dual mappings to the homomorphisms $T_n:\ S^{(\nu )}\to S^{(n)}$.
Using the duality theorem [16], we find that
\begin{equation}
S^*=\varinjlim K_n/O_n
\end{equation}
where the direct limit is taken with respect to the imbeddings.
The right-hand side of (4) equals $K/O$. $\qed $

\vspace{1pc}
Note that the above isomorphism is given by the same formula (3)
where this time $a_\psi $ is an arbitrary representative of a
coset from $K/O$.

Let $\varphi _a(x)=\chi (<a,x>),\ x\in \K $, where $a\in K,\
\Vert a\Vert >1$ or $a=0$. It is known [1] that
\begin{equation}
(D^\alpha \varphi _a)(x)=\Vert a\Vert ^\alpha \varphi _a(x)
\end{equation}
for $\mu $-almost all $x\in \K$. Note that the set of values of
the function $a\mapsto \Vert a\Vert ^\alpha$ with $\Vert a\Vert
>1$ coincides with the set
\begin{equation}
\left\{ q_n^{\alpha N/m_n},\ n,N=1,2,\ldots \right\} =
\left\{ q_1^{\alpha N/e_n},\ n,N=1,2,\ldots \right\}
\end{equation}
where $e_n$ is the ramification index of the extension $K_n/k$.
Denote its residue degree by $f_n$. It is well known that the
sequences $\{ f_n\}$, $\{ e_n\}$ are non-decreasing and
$e_nf_n=m_n$.

\vspace{2pc}
\noindent THEOREM 1\quad
Let $A\subset K$ be a complete system of representatives of
cosets from $K/O$. Then $\{ \varphi _a\} _{a\in A}$ is the
orthonormal eigenbasis for the operator $D^\alpha $ in
$L_2(S,d\mu )$. As a set, its spectrum equals the set (6)
complemented with the point $\lambda =0$. Each non-zero
eigenvalue of $D^\alpha $ has an infinite multiplicity.
The point $\lambda =0$ is an accumulation point for eigenvalues
if and only if $e_n\to \infty $.

\vspace{1pc}
\noindent Proof: \quad
The first statement follows from (5) and Lemma 1. The assertion
about the accumulation at zero is obvious from (6).

Let us construct $A$ as the union of an increasing family $\{
A_n\}$ of complete systems of representatives of
cosets from $K_n/O_n$. Each $A_n$ consists of elements of the
form $a=\pi _n^{-N}\left( \sigma _1+\sigma _2\pi _n+\cdots
+\sigma _{N-1}\pi _n^{N-1}\right) $, $N\ge 1$, where $\sigma
_1,\ldots ,\sigma _{N-1}$ belong to the set $F_n$ of
representatives of the residue field corresponding to the field
$K_n$, $\sigma _1\ne 0$, $\pi _n$ is a prime element of $K_n$.
We have $|a|_n=q_n^N$ so that $\Vert
a\Vert =q_1^{Nf_n/m_n}=q_1^{N/e_n}$. Meanwhile $\mbox{card
}F_n=q_n=q_1^{f_n}$.

If $e_n\to \infty $ then the same value of $\Vert a\Vert $
corresponds to elements from infinitely many different sets $A_n$
with different values of $N$ ($e_n$ is a multiple of $e_{n-1}$
due to the chain rule for the ramification indices; see [19]).
If the sequence $\{ e_n\}$ is
bounded, then it must stabilize, and we obtain the same value
of $\Vert a\Vert $ for elements from infinitely many sets $A_n$
with possibly the same $N$. However, in this case $f_n\to \infty
$, the number of such elements (for a fixed $N$) is $Nq_1^{f_n}-
1$ ($\to \infty $ for $n\to \infty $). In both cases we see that
all the non-zero eigenvalues have an infinite multiplicity. $\qed
$

\vspace{1pc}
Note that the cases where $e_n\to \infty $ or $e_n\le
\mbox{const}$ both appear in important examples of infinite
extensions. Let $K$ be the maximal unramified extension of $k$.
Then one may take for $K_n$ the unramified extension of $k$ of
the degree $n!$, $n=2,3,\ldots $. Here $e_n=1$ for all $n$. On
the other hand, if $K$ is the maximal abelian extension of
$k=Q_p$ then a possible choice of $K_n$ is the cyclotomic
extension $C_{n!}=Q_p(W_{n!})$ where $W_l$ is the set of all
roots of 1 of the degree $l$ (see [20,21]). Writing
$n!=n'p^{l_n}$, $(n',p)=1$, we see that $e_n=(p-
1)p^{l_n-1}\to \infty $ as $n\to \infty$.

\vspace{3pc}
\noindent {\bf 3\quad HYPERSINGULAR INTEGRAL REPRESENTATION}

\vspace{2pc}
The main aim of this section is the following result.

\vspace{1pc}
\noindent THEOREM 2\quad
There exists such a measure $\Pi $ on $\mathcal B(\K \setminus \{
0\} )$ finite outside any neighbourhood of the origin that
$D^\alpha $ has the representation (2) on all functions $f\in
\mathcal E(\K )$.

\vspace{1pc}
\noindent
In the course of proof we shall also obtain some new information
about the Markov process $X(t)$ generated by the operator
$-D^\alpha $. We assume below that $X(0)=0$.

The projective limit topology on $\K$ coincides with the one
given by the shift-invariant metric
$$
r(x,y)=\sum \limits _{n=1}^\infty 2^{-n}\frac{\Vert x-y\Vert
_n}{1+\Vert x-y\Vert _n},\quad x,y\in \K .
$$
It is known [22] that the main notions and results
regarding stochastic processes with independent increments carry
over to the case of a general topological group with a
shift-invariant metric.

Let $\nu (t,\Gamma )$ be a Poisson random measure corresponding
to the process $X(t)$. Here $\Gamma \in \mathcal B_0=\bigcup
\limits _{\gamma >0}\mathcal B_\gamma $,
$$
\mathcal B_\gamma =\{ \Gamma \subset \mathcal B(\K ),\
\mbox{dist}(\Gamma ,0)\ge \gamma \}.
$$
For any $\Gamma \in \mathcal B_0$ we can define a stochastically
continuous process $X_\Gamma (t)$, the sum of all jumps of the
process $X(\tau )$ for $\tau \in [0,t)$ belonging to $\Gamma $.
In a standard way [23] we find that
\begin{equation}
E\chi (<\lambda ,X_\Gamma (t)>)=\exp \left\{ \int \limits _\Gamma
(\chi (<\lambda ,x>)-1)\Pi (t,dx)\right\} ,\quad \lambda \in K,
\end{equation}
where $E$ denotes the expectation, $\Pi (t,\cdot )=E\nu (t,\cdot
)$.

Let $\lambda \in K_n$. Consider the set
$$
V_{\delta ,n}=\{ x\in \K :\ \Vert T_n(x)\Vert \ge \delta \}
,\quad 0<\delta \le 1.
$$
Let us look at the equality (7) with $\Gamma =V_{\delta ,n},\
\delta \le \Vert \lambda \Vert ^{-1}$. If $x\in V_{\delta ,n}$
then $r(x,0)\ge 2^{-n}\frac{\Vert x\Vert _n}{1+\Vert x\Vert
_n}\ge 2^{-n-1}\delta $, whence $V_{\delta ,n}\in \mathcal B_0$.
If $\delta \le \Vert \lambda \Vert ^{-1}$, $x\notin V_{\delta
,n}$, then
$$
\Vert <\lambda ,x>\Vert =\Vert T(\lambda T_n(x))\Vert \le 1
$$
which implies that the integral in the right-hand side of (7)
coincides with the one taken over $\K$.

On the other hand, almost surely
\begin{equation}
\Vert X(t)-X_{V_{\delta ,n}}(t)\Vert _n<\delta .
\end{equation}
Indeed, let $t_0$ be the first exit time of the process
$X(t)-X_{V_{\delta ,n}}(t)$ from the set $\K \setminus V_{\delta
,n}$. Suppose that $t_0<\infty $ with a non-zero probability.
Then
\begin{equation}
\Vert X(t_0-0)-X_{V_{\delta ,n}}(t_0-0)\Vert _n<\delta ,
\end{equation}
\begin{equation}
\Vert X(t_0+0)-X_{V_{\delta ,n}}(t_0+0)\Vert _n\ge \delta .
\end{equation}

If  $\Vert X(t_0+0)-X(t_0-0)\Vert _n<\delta $ then
$X_{V_{\delta ,n}}(t_0+0)=X_{V_{\delta ,n}}(t_0-0)$, so that
\begin{equation}
\Vert [X(t_0+0)-X_{V_{\delta ,n}}(t_0+0)]-
[X(t_0-0)-X_{V_{\delta ,n}}(t_0-0)]\Vert _n<\delta .
\end{equation}
If, on the contrary, $\Vert X(t_0+0)-X(t_0-0)\Vert _n\ge \delta
$, then
$$
X_{V_{\delta ,n}}(t_0+0)-X_{V_{\delta ,n}}(t_0-0)=X(t_0+0)-X(t_0-
0),
$$
the expression in the left-hand side of (11) equals zero, and the
inequality (11) is valid too. In both cases it contradicts (9),
(10). Thus $t_0=\infty $ almost surely, and the inequality (8)
has been proved. We have come to the following formula of
L\'evy-Khinchin type.

\vspace{1pc}
\noindent LEMMA 2\quad
For any $\lambda \in K,\ t>0$
\begin{equation}
E\chi (<\lambda ,X(t)>)=\exp \left\{ \int \limits _{\K}[\chi
(<\lambda ,x>)-1]\Pi (t,dx)\right\} .
\end{equation}

\vspace{1pc}
\noindent REMARK \quad
Lemma 2 can serve as a base for developing the theory of
stochastic integrals and stochastic differential equations over
$\K$. In fact, the techniques and results of [24]
carry over to this case virtually unchanged.

\vspace{1pc}
Both sides of (12) can be calculated explicitly if we use the
heat measure $\pi (t,dx)$ (see [1]). We have
$$
E\chi (<\lambda ,X(t)>)=\int \limits _{\K}\chi (<\lambda ,x>)\pi
(t,dx)=\rho _\alpha (\Vert \lambda \Vert ,t)
$$
where
\begin{equation}
\rho _\alpha (s,t)=\left\{ \begin{array}{rl}
e^{-ts^\alpha }, & \mbox{if } s>1 \\
1, & \mbox{if } s\le 1,
\end{array} \right.
\end{equation}
$s\ge 0,\ t>0$. It follows from the definitions that $\Pi
(t,\cdot )$ is symmetric with respect to the reflection $x\mapsto
-x$. Therefore
\begin{equation}
\int \limits _{\K}[\chi (<\lambda ,x>)-1]\Pi (t,dx)=
\left\{ \begin{array}{rl}
-t\Vert \lambda \Vert ^\alpha , & \mbox{if } \Vert \lambda \Vert >1 \\
0, & \mbox{if } \Vert \lambda \Vert \le 1,
\end{array} \right.
\end{equation}

\vspace{1pc}
\noindent LEMMA 3\quad
Let $M_n$ be a compact subset of $K_n\setminus \{ 0\}$, $M=T_n^{-
1}(M_n)$. Then
\begin{equation}
\Pi (t,M)=-t\int \limits _{\eta \in K_n:\ \Vert \eta \Vert >1}
\Vert \eta \Vert ^\alpha w_n(\eta )\,d\eta
\end{equation}
where $w_n(\eta )$ is the inverse Fourier transform of the
function $y\mapsto q_n^{-d_n/2}\omega _{M_n}^{(n)}(m_ny)$,
$\omega _{M_n}^{(n)}$ is the indicator of the set $M_n$ in $K_n$,
and $dx$ is the normalized additive Haar measure on $K_n$.

\vspace{1pc}
\noindent Proof: \quad
Let $\omega _M$ be the indicator of the set $M$ in $\K$. Then
$\omega _M(x)=\omega _{M_n}^{(n)}(T_n(x))$, $x\in \K$,
$$
\omega _{M_n}^{(n)}(\xi )=\int \limits _{K_n}\chi (\xi \eta
)w_n(\eta )\,d\eta ,\quad \xi \in K_n.
$$
Since
$$
\int \limits _{K_n}w_n(\eta )\,d\eta =\omega _{M_n}^{(n)}(0)=0,
$$
we get
$$
\omega _M(x)=\int \limits _{K_n}[\chi (\eta T_n(x))-1]w_n(\eta
)\,d\eta .
$$
Integrating with respect to $\Pi (t,dx)$ and using (14) we come
to (15). $\qed $

\vspace{1pc}
It follows from Lemma 3 that $\Pi (t,dx)=t\Pi (1,dx)$. We shall write
$\Pi (dx)$ instead of $\Pi (1,dx)$.

\vspace{1pc}
\noindent Proof of Theorem 2: \quad
Let $f(x)=\varphi (T_n(x)),\ x\in \K$, where $\varphi $ is
locally constant and in addition supp $\varphi $ is compact,
$0\notin \mbox{supp }\varphi $. It follows from Lemma 3 that
\begin{equation}
\int \limits _{\K}f(x)\Pi (dx)=-\int \limits _{K_n}\Delta ^\alpha
(\eta )\psi (\eta )\,d\eta
\end{equation}
where $\psi $ is the inverse Fourier transform of the function
\linebreak $y\mapsto q_n^{-d_n/2}\varphi (m_ny)$.

The right-hand side of (16) is an entire function with respect to
$\alpha $. Assuming temporarily Re $\alpha <-1$, we can use the
Plancherel formula with subsequent analytic continuation (see
[1]). As a result we find that for $\alpha >0$
\begin{multline}
\int \limits _{\K}f(x)\Pi (dx) =
-q_n^{d_n\alpha /m_n}\frac{1-q_n^{\alpha /m_n}}{1-q_n^{-1-
\alpha /m_n}}\int \limits_{x\in K_n:\  |x|_n
\le q_n^{d_n}}\biggl[ |x|_n^{-1-\alpha /m_n}\biggr. \\ +\biggl.
\frac{1-q_n^{-1}}{q_n^{\alpha /m_n}-1}q_n^{-d_n(1+\alpha /m_n)}
\biggr] \varphi (m_nx)\,dx.
\end{multline}

An obvious approximation argument shows that (17) is valid for
any $f\in \mathcal E(\K )$. Comparing (17) with (1) we obtain
(2). $\qed $

\vspace{3pc}
\noindent {\bf 4\quad HEAT MEASURE}

\vspace{2pc}
Recall that the heat measure $\pi (t,dx)$ corresponding to the
operator $-D^\alpha $ is defined by the formula
$$
\int \limits _{\K}\chi (<\lambda ,x>)\pi (t,dx)=\rho _\alpha (\|
\lambda \|,t),\quad \lambda \in K,\ t>0,
$$
where $\rho _\alpha$ is given by (13).

\vspace{1pc}
\noindent THEOREM 3\quad
For each $t>0$ the measure $\pi (t,\cdot )$ is not absolutely
continuous with respect to $\mu $.

\vspace{1pc}
\noindent Proof: \quad
Let us fix $N\ge 1$ and consider the set
$$
M=\left\{ x\in \overline{K}:\ \ \|T_n(x)\|\le q_n^{d_n/m_n-N/f_n}\|m_n\|,\
n=1,2,\ldots \right\}
$$
We shall show that $\mu (M)=0$ whereas $\pi (t,M)\ne 0$.

Denote
$$
M_n=\left\{ x\in \overline{K}:\ \ \|T_n(x)\|\le q_n^{d_n/m_n-N/f_n}\|m_n\|
\right\} ,\quad n=1,2,\ldots .
$$
It is clear that $M=\bigcap \limits _{n=1}^\infty M_n$. Repeating
the arguments from the proof of Lemma 1, we see that $M_\nu
\subset M_n$ if $\nu >n$. Thus
$$
\pi (t,M)=\lim \limits _{n\to \infty }\pi (t,M_n),\quad \mu
(M)=\lim \limits _{n\to \infty }\mu (M_n).
$$

It follows from the integration formula for cylindrical functions
[1] that
\begin{multline}
\mu (M_n)=q_n^{-d_n}\|m_n\|^{-m_n}\int \limits _{z\in K_n:\ \|z\|\le
q_n^{d_n/m_n-N/f_n}\|m_n\|}dz\\
= q_n^{-d_n}|m_n|_n^{-1}\int \limits _{z\in K_n:\ |z|_n\le
q_n^{d_n-Ne_n}|m_n|_n}dz=q_n^{-Ne_n}=q_1^{-Nf_ne_n},
\end{multline}
so that $\mu (M_n)=q_1^{-Nm_n}\to 0$ for $n\to \infty $. Thus
$\mu (M)=0$.

In a similar way (see [1])
\begin{equation}
\pi (t,M_n)=|m_n|^{-1}_n\int \limits _{z\in K_n:\ |z|_n\le
q_n^{d_n-Ne_n}|m_n|_n}\Gamma _\alpha ^{(n)}\left( m_n^{-1}z,t\right)
\,dz
\end{equation}
where $\Gamma _\alpha ^{(n)}$ is a fundamental solution of the
Cauchy problem for the equation over $K_n$ of the form
$\frac{\partial u}{\partial t}+\partial _n^\alpha u=0$. Here
$\partial _n^\alpha$ is a pseudo-differential operator over $K_n$
with the symbol $\Delta ^\alpha (\xi )$. It is clear that
$$
\Gamma _\alpha ^{(n)}(\zeta ,t)=q_n^{-d_n/2}\widetilde{\rho
}_\alpha (\zeta ,t)
$$
where tilde means the local field Fourier transform:
$$
\widetilde{u}(\zeta )=q_n^{-d_n/2}\int _{K_n}\chi \circ \mbox{Tr}_
{K_n/k}(z\zeta )u(z)\,dz\,,
$$
for a complex-valued function $u$ over $K_n$ (sufficient
conditions for the existence of $\widetilde{u}$ and the validity
of the inversion formula
$$
u(z)=q_n^{-d_n/2}\int _{K_n}\chi \circ \mbox{Tr}_{K_n/k}(-z\zeta
)\widetilde{u}(\zeta )\,d\zeta \,,
$$
are well known).

Using the Plancherel formula we can rewrite (19) in the form
$$
\pi (t,M_n)=\int \limits _{K_n}\widetilde{\Gamma }^{(n)}_\alpha
(x,t)\widetilde{\beta }_n(x)\,dx
$$
where
$$
\beta _n(\zeta )=\left\{ \begin{array}{rl}
1, & \mbox{if } |\zeta |_n\le q_n^{d_n-Ne_n} \\
0, & \mbox{if } |\zeta |_n> q_n^{d_n-Ne_n}.
\end{array} \right.
$$

We have $\widetilde{\Gamma }^{(n)}_\alpha =q_n^{-d_n/2}\rho
_\alpha $,
$$
\widetilde{\beta }_n(x)=\left\{ \begin{array}{rl}
q_n^{\frac{d_n}2-Ne_n}, & \mbox{if } |x|_n\le q_n^{Ne_n} \\
0, & \mbox{if } |\zeta |_n> q_n^{Ne_n}
\end{array} \right.
$$
(see e.g. [11]), so that
\begin{multline*}
\pi (t,M_n)=q_n^{-Ne_n}\int \limits _{|x|_n\le q_n^{Ne_n}}
\rho _\alpha (\|x\|,t)\,dx\\ \ge
q_n^{-Ne_n} \int \limits _{|x|_n=q_n^{Ne_n}}\rho _\alpha \left(
|x|_n^{1/m_n},t\right) \,dx
=\left( 1-q_n^{-1}\right) \exp \left( -tq_n^{\alpha N/f_n}\right)\\ \ge
\left( 1-q_1^{-1}\right) \exp \left( -tq_1^{\alpha N}\right) .
\end{multline*}
Hence, $\pi (t,M)>0$. $\qed$

\vspace{2pc}
\noindent ACKNOWLEDGEMENT\quad
This work was supported in part by the Ukrainian Fund for
Fundamental Research (Grant 1.4/62).

\vspace{3pc}
\noindent {\bf REFERENCES}

\vspace{2pc}\noindent
\noindent \ 1. AN Kochubei. Analysis and probability over
infinite extensions of a local field. Potential Anal. (to
appear).

\noindent \ 2. SN Evans. Local field Gaussian measures. In:
E Cinlar, KL Chung, RK Getoor, eds.  Seminar on Stochastic
Processes 1988. Boston: Birkh\"auser, 1989, pp. 121-160.

\noindent \ 3. VS Vladimirov. Generalized functions over the field
of p-adic numbers. Russian Math. Surveys 43, No. 5: 19-64, 1988.

\noindent \ 4. VS Vladimirov, IV Volovich, EI Zelenov.
p-Adic Analysis and Mathematical Physics. Singapore: World Scientific,
1994.

\noindent \ 5. RS Ismagilov. On the spectrum of the self-adjoint
operator in $L_2(K)$ where $K$ is a local field; an analog of the
Feynman-Kac formula. Theor. Math. Phys. 89: 1024-1028, 1991.

\noindent \ 6. S Haran. Riesz potentials and explicit sums
in arithmetic. Invent. Math. 101: 697-703, 1990.

\noindent \ 7. S Haran. Analytic potential theory over the p-
adics. Ann. Inst. Fourier 43: 905-944, 1993.

\noindent \ 8. AN Kochubei. Schr\"odinger-type operator
over p-adic number field. Theor. Math. Phys. 86: 221-228, 1991.

\noindent \ 9. AN Kochubei. Parabolic equations over the field of p-adic
numbers. Math. USSR Izvestiya 39: 1263-1280, 1992.

\noindent 10. AN Kochubei. The differentiation operator on
subsets of the field of p-adic numbers. Russ. Acad. Sci.
Izv. Math. 41: 289-305, 1993.

\noindent 11. AN Kochubei. Gaussian integrals and spectral theory
over a local field. Russ. Acad. Sci. Izv. Math. 45: 495-503,
1995.

\noindent 12. AN Kochubei. Heat equation in a p-adic ball.
Methods of Funct. Anal. and Topology 2, No. 3-4: 53-58, 1996.

\noindent 13. AD Blair. Adelic path space integrals. Rev.
Math. Phys. 7: 21-50, 1995.

\noindent 14. VS Varadarajan. Path integrals for a class
of p-adic Schr\"odinger equations. Lett. Math. Phys. 39: 97-106,
1997.

\noindent 15. AD Bendikov. Symmetric stable semigroups on
the infinite-dimen\-sio\-nal torus. Expositiones Math. 13: 39-79,
1995.

\noindent 16. S Kaplan. Extensions of the Pontrjagin
duality. II: Direct and inverse sequences. Duke Math. J. 17:
419-435, 1950.

\noindent 17. JWS Cassels, A Fr\"ohlich, eds. Algebraic Number Theory.
London and New York: Academic Press, 1967.

\noindent 18. A Weil. Basic Number Theory. Berlin: Springer,
1967.

\noindent 19. IB Fesenko, SV Vostokov. Local Fields and Their Extensions.
Providence: American Mathematical Society, 1993.

\noindent 20. K Iwasawa. Local Class Field Theory. Tokyo: Iwanami
Shoten, 1980 (in Japanese; Russian translation, Moscow: Mir, 1983).

\noindent 21. JP Serre. Local Fields. New York: Springer, 1979.

\noindent 22. AV Skorohod. Random Processes with
Independent Increments. Dordrecht: Kluwer, 1991.

\noindent 23. II Gikhman, AV Skorohod. Theory of
stochastic processes II. Berlin: Springer, 1975.

\noindent 24. AN Kochubei. Stochastic integrals and
stochastic differential equations over the field of p-adic
numbers. Potential Anal. 6: 105-125, 1997.

\end{document}